\documentclass[12pt,reqno]{amsart}
\usepackage{mathrsfs}
\usepackage{amsgen}
\usepackage{amscd}
\usepackage{amsmath}
\usepackage{latexsym}
\usepackage{amsfonts}
\usepackage{amssymb}
\usepackage{amsthm}
\usepackage{graphicx}
\usepackage{color}
\usepackage{amsthm,amssymb}
\usepackage{graphicx}
\usepackage{enumerate}
\usepackage{amssymb,amsmath,graphicx,amsfonts,euscript}
\usepackage{color}
\usepackage{mathrsfs}
\usepackage{empheq}
\usepackage{cases}
\usepackage{cite}
\usepackage{ulem}
\usepackage{cancel}

\usepackage{hyperref}

\usepackage[margin=3cm, a4paper]{geometry}

\def\H{{\cal H}}

\def\R{\mathbb{R}}

\def\T{\mathbb{T}}
\def\C{\mathbb{C}}

\def\H2{H^2(\R^N)}
\def\L2{L^2(\R^N)}
\def\to{\rightarrow}

\newcommand{\dx}{\,\mathrm{d}x}

\newcommand{\dt}{\,\mathrm{d}t}

\def\H{{\cal H}}

\def\H1{H^1(\R)}

 \newcommand{\Del}[1]{}

\numberwithin{equation}{section}

\newtheorem{thm}{Theorem}[section]

\newtheorem{lem}[thm]{Lemma}

\theoremstyle{remark}
\newtheorem{remark}[thm]{Remark}
\newtheorem*{exam*}{Examples}

\begin{document}

\setcounter{page}{1}

\title[Scattering for gDNLS]{Optimal small data Scattering for the generalized derivative nonlinear schr\"odinger equations}

\author{Ruobing Bai}
\address{Center for Applied Mathematics\\
Tianjin University\\
Tianjin 300072, China}
\email{baimaths@hotmail.com}
\thanks{RB and YW are partially supported by NSFC 11771325 and 11571118. YW is also supported by the National Youth Topnotch Talent Support Program in China. JX is supported by the Research Council of Norway (No. 250070)}

%
%\author{Zihua Guo}
%\address{School of Mathematical Sciences, Monash University, VIC 3800, Australia}
%\email{zihua.guo@monash.edu}
%\thanks{Z. Guo is supported by ARC DP170101060.}
%

\author{Yifei Wu*}
\address{Center for Applied Mathematics, Tianjin University,
Tianjin 300072, China}
\email{yerfmath@gmail.com}
\thanks{* Corresponding author}

\author{Jun Xue}
\address{Department of Mathematical Sciences\\
Norwegian University of Science and Technology\\
Trondheim 7491, Norway}
\email{jxuemath@hotmail.com}
\thanks{\href{https://doi.org/10.1016/j.jde.2020.05.001}{Journal of Differential Equations, Volume 269, Issue 9, Pages 6422--6447}. This paper is the final version.}


%\date{\today}


\maketitle

\begin{abstract}\noindent
In this work, we consider the following generalized derivative nonlinear Schr\"odinger equation
\begin{align*}
 i\partial_t u+\partial_{xx} u +i |u|^{2\sigma}\partial_x u=0, \quad (t,x)\in \R\times \R.
\end{align*}
We prove that when $\sigma\ge 2$, the solution is global and scattering when the initial data is small in $H^s(\R)$,  $\frac 12\leq s\leq1$.
Moreover, we show that when $0<\sigma<2$, there exist a class of solitary wave solutions $\{\phi_c\}$ satisfying
$$
\|\phi_c\|_{H^1(\R)}\to 0,
$$
when $c$ tends to some endpoint,  which is against the small data scattering statement. Therefore, in this model, the exponent $\sigma\ge2$ is optimal for small data scattering. We remark that this exponent is larger than the short range exponent and the Strauss exponent.
\end{abstract}

\section{Introduction}
\vskip 0.2cm
In this paper, we consider the small data scattering of the Cauchy problem for the following generalized derivative nonlinear Schr\"odinger equation (gDNLS)
\begin{align}\label{GDNLS}
 \left\{ \begin{aligned}
 &i\partial_t u+\partial_x^2 u +i |u|^{2\sigma}\partial_x u=0, \qquad (t,x)\in \R\times \R,\\
&u(x,0)=\varphi(x). \\
   \end{aligned}
  \right.
\end{align}
Here $\sigma > 0$, $u:\R \to \C$ is an unknown function.

The generalized derivative nonlinear Schr\"odinger equation describes the physical phenomenon of Alfv\'en waves with small but finite amplitude propagating along the magnetic field in cold plasmas (see for example \cite{MioOginoMinamiTakeda-JPS-76}).

The class of solutions to equation (\ref{GDNLS}) is invariant under the scaling
\begin{equation}\label{eqs:scaling-p}
u(t,x)\to u_\lambda(t,x) = \lambda^{\frac1{2\sigma}} u(\lambda^2 t, \lambda x) \ \ {\rm for}\ \ \lambda>0,
\end{equation}
which maps the initial data as
\begin{eqnarray}
u(0)\to u_{\lambda}(0):=\lambda^{\frac1{2\sigma}} \varphi(\lambda x) \ \ {\rm for}\ \ \lambda>0.\nonumber
\end{eqnarray}
Denote
$$
s_c=\frac 12-\frac1{2\sigma},
$$
then the scaling  leaves  $\dot{H}^{s_{c}}$ norm invariant, that is,
\begin{eqnarray*}
\|u\|_{\dot H^{s_{c}}_x}=\|u_{\lambda}\|_{\dot H^{s_{c}}_x}.
\end{eqnarray*}

When $\sigma=1$, take a suitable gauge transformation
$$u(t,x)\to u(t,x) \exp\Big(-\frac i2 \int_{-\infty}^x |u(t,y)|^2\,\mathrm dy \Big),$$
then the equation in \eqref{GDNLS}
is transformed into the standard derivative nonlinear Schr\"odinger equation (DNLS)
\begin{align}\label{DNLS}
i\partial_t u+\partial_x^2 u +i\partial_x(|u|^2u) =0.
\end{align}
For $\sigma \neq 1$, \eqref{GDNLS} is regarded as a generalization of \eqref{DNLS}.
The well-posedness theory and the long time behavior of the solution for the equation \eqref{DNLS} have been widely considered by many researchers. For the local well-posedness result, Hayashi and Ozawa \cite{HaOza-PhyD-92,HaOza-SJMA-94} proved that equation \eqref{DNLS} is locally well-posed in the Sobolev space $H^1(\R)$ (see also the previous works \cite{GuoTan-PRSE-91, TuFu-80-DNLS}).
Later, the best result regarding local well-posedness was obtained in $H^{\frac 12}(\mathbb R)$ by Takaoka in \cite{Ta-99-DNLS-LWP}.
Very recently, Mosincat and Yoon \cite{MoYo} proved the unconditional well-posedness in $H^s(\R), s>\frac12$ (see also Dan, Li and Ning \cite{DanNiLi} for the previous work in $H^s(\R), s>\frac23$).
Some other results on local well-posedness can be found in \cite{Gr-05,Grhe-95,Herr,Ta-16-DNLS-LWP}. Moreover, Biagioni and Linares \cite{BiLi-01-Illposed-DNLS-BO} proved the equation \eqref{DNLS} is ill-posed in $H^s(\mathbb R), s<\frac 12$.
With regard to the theory of global well-posedness, Hayashi and Ozawa \cite{HaOza-PhyD-92} proved that it is globally well-posed in $H^1(\R)$ under the condition that the initial data satisfies $\|u_0\|_{L^2}< \sqrt {2\pi}$.
Wu \cite{Wu-APDE-13,Wu-APDE-15} showed that it is globally well-posed in $H^1(\R)$ under the condition $\|u_0\|_{L^2}< 2\sqrt \pi$.
Guo and Wu \cite{GuoWu-DCDS-17} later proved that it is globally well-posed in $H^{\frac 12}(\R)$ under the same condition of initial data (see also \cite{CollianderTao-SJMA-01, CollianderTao-SJMA-02, Miao-Wu-Xu:2011:DNLS,Ta-01-DNLS-GWP} for the previous results on the low regularity). The same results also hold in the periodic case (see Mosincat and Oh \cite{MoOh} in $H^1(\T)$, and Mosincat \cite{Mosincat} in $H^\frac12(\T)$). More recently, Jenkins, Liu, Perry and Sulem \cite{JeLiPeSu} proved that the Cauchy problem \eqref{DNLS} is globally well-posed in the weighted Sobolev space $H^{2,2}(\R)$.

The equation in \eqref{GDNLS} in the case of $\sigma \not =1$ also attracts a lot of researchers in recent years. Firstly, for the local well-posedness result, when $0<\sigma<\frac 12$, Linares, Ponce and Santos
\cite{LinaresPonceSantos-PRE,LinaresPonceSantos-PRE-II}  proved the local well-posedness for a class of data of arbitrary size in an appropriate weighted
Sobolev space. When $\frac 12 \leq \sigma< 1$, Hayashi and Ozawa \cite{HaOza-JDE-16}
proved that (gDNLS) is locally well-posed in $H^2(\R)$, and Santos \cite{Santos-JDE-2015}  showed the local well-posedness in a weighted sobolev space. When $\sigma >1$, Hayashi and Ozawa \cite{HaOza-JDE-16} proved that (gDNLS) is locally well-posed in energy space
$H^1(\R)$. Hao \cite{Hao-CPAA-07} proved that it is locally well-posed in $H^\frac 12(\R)$, when $\sigma \ge \frac 52$. Santos \cite{Santos-JDE-2015} proved that it is locally well-posed in $H^\frac 12(\R)$ with small initial data when $\sigma >1$.
Secondly, compared with the local well-posedness, there are only  few results of global well-posedness. When $0<\sigma <1$, Hayashi and Ozawa \cite{HaOza-JDE-16} showed the global existence without uniqueness of (gDNLS) in $H^1(\R)$. When $\sigma > 1$, Fukaya, Hayashi and Inui \cite{FukayaHayashiInui-APDE-17} gave a sufficient condition of initial data for global well-posedness in $H^1(\R)$. Some other results related to the stability theory and inverse scattering theory can be found in \cite{ChSiSu-DNLS, CoOh-06-DNLS,  Fu-16-DNLS, GuWu95, guo-DNLS,GuoNingWu-Pre, JeLiPeSu-1, KwonWu-JAM-18, Stefan-W-15-MultiSoliton-DNLS, LiBing-NingCui-PRE, LiSiSu1, LiPeSu1, LiPeSu3, LiSiSu, Miao-Tang-Xu:2016:DNLS, Miao-Tang-Xu:2017:DNLS, PeSh-DNLS, PeSh-DNLS-2, TaXu-17-DNLS-Stability} and the references therein.

The equation that we investigate in \eqref{GDNLS} in this paper can be also treated as the form of
\begin{align}
i\partial_tu+\Delta u=P(u,\bar u,\partial_xu,\partial_x\bar u).\label{14:27}
\end{align}
The well-posedness theory of equation \eqref{14:27} has been studied by many researchers. Here we only introduce some works and readers can seek further literatures from related references, when $P$ is a polynomial of the form $P(z)=\sum_{d\leq|\alpha|\leq l }C_\alpha z^\alpha$ and $l, d$ are integers with $l\ge d$. For general cases with $d\ge 3$, Kenig, Ponce and Vega \cite{KenigPonceVega-AIPANL-93} showed that the equation \eqref{14:27} is locally well-posed with small initial data in $H^{\frac 72}(\R)$. Some further results have been acquired when $P$ is only composed of $\bar u$ and $\partial_x \bar u$ under some suitable assumption. Gr\"unrock \cite{Grun-PRE} proved that the equation \eqref{14:27} is locally well-posed for $s>\frac 12-\frac 1{d-1}$ when $P=\partial_x(\bar u^d)$ and $s>\frac 32-\frac 1{d-1}$ when $P=(\partial_x\bar u)^d$ respectively. Hirayama \cite{Hira-FE-15} later extended Gr\"unrock's results to the small data global well-posedness for
$s\ge \frac 12-\frac 1{d-1}$ when $P=\partial_x(\bar u^d)$. Recently, Pornnopparath \cite{Pornnopparath-PRE} proved that when each term in $P$ contains only one derivative, the equation \eqref{14:27} is locally well-posed in $H^\frac 12(\R)$, and when a term in $P$ has more than one derivative, the equation \eqref{14:27} is locally well-posed in $H^\frac 32(\R)$.
Moreover, Pornnopparath also proved that when $d\ge 5$, \eqref{14:27} is almost globally well-posed in $H^s(\R)$ when $P$ has only one
derivative and $s> \frac 12$, or when $P$ has more than one derivative and $s>\frac 32$.
For higher dimension and more related theories, see \cite{Bejenaru1-IMRP-06, Bejenaru1-IMRP-08, Chihara-MA-99, KenigPonceVega-AIPANL-93, WangY-CPDE-11} and the references therein.

All of the results above are related to the theories of local and global well-posedness. To our knowledge, there is no scattering result yet to (gDNLS). The related result on modified scattering can be found in \cite{GuHaLinNa-DNLS, HaOza-MA-94} and the references therein.

One of the motivation to prove scattering is that we believe the small data scattering result of the present paper is significative to our further study.
 In order to consider the long-time behavior of the solution
to (gDNLS), the small data scattering theory is initially needed in some situation, for example, long-time perturbation theory when we use the
concentration-compactness argument.

Moreover, it was known that when $\sigma=1$, there exist solitary wave solutions which can be arbitrarily close to zero. This implies that the small data scattering is not true when $\sigma=1$. So one may wonder the optimal value of $\sigma$ such that the scattering statement holds when the initial data is small enough in some Sobolev space. This is another motivation in the present paper.

For semilinear Schr\"odinger equation, there are two important exponents named short range exponent and the Strauss exponent. When the nonlinear power is larger than the short range exponent $3$ ($1+\frac2d$ for general dimensions), one has the global well-posedenss and the existence of the wave operator  for small data (see for examples \cite{CaWe-1992-CMP, GiOzVe-1994-ANPoincare, Nakanishi-SIAM-2001}); when the nonlinear power is larger than the Strauss exponent $\frac{\sqrt{17}+3}{2}\approx 3.56$ ($\frac{\sqrt{d^2+12d+4}+d+2}2$ for general dimensions), one has the scattering for small data (see \cite{Strauss-JFA-1981}). According to these, especially because of the short range exponent, one may ask whether $\sigma=1$ is the optimal exponent for scattering. However, there is no such general result for non-semilinear Schr\"odinger equation (related results see \cite{Cohn, Delort, GeMaSh, HaMiNa} and the references therein). In fact, it is of much model dependence when the nonlinearity contains derivatives. In the present paper, as what we will show in the following, the situation for the nonlinear Schr\"odinger equation with derivatives is of much difference, compared with the semilinear Schr\"odinger equations, and the models mentioned in the references above, the optimal exponent for scattering is 5, which is much larger than the short range exponent and the Strauss exponent.

For all $\sigma>0$, the equation \eqref{GDNLS} has a two-parameter family of solitary waves,
$$
u_{\omega,c}(t)=e^{i\omega t}\phi_{\omega, c}(x-ct),
$$
where the parameters $c^2<4\omega$, and $\phi_{\omega, c}$ is the solution of the form
\begin{align}
\phi_{\omega,c}(x)=\varphi_{\omega,c}(x)\exp{\big\{\frac c2 i x-\frac{i}{2\sigma+2}\int_{-\infty}^{x}\varphi^{2\sigma}_{\omega,c}(y)dy\big\}},\label{phi}
\end{align}
with
\begin{align*}
\varphi_{\omega,c}(x)=\Big\{\frac{(\sigma+1)(4\omega-c^2)}{2\sqrt\omega \cosh(\sigma\sqrt{4\omega-c^2}\,x)-c}\Big\}^{\frac{1}{2\sigma}}.
\end{align*}

Firstly, we show that there exists a sequence of solitary waves which are arbitrary small in $H^1(\R)$ when $0<\sigma<2$; while all of the solitary waves are away from the origin when $\sigma\ge 2$.  Let
$\Omega=\{(\omega,c): c\in(-2\sqrt \omega, 2\sqrt \omega)\}$.
\begin{thm}\label{thm:main1}
Let $\phi_{\omega,c}$ be defined in \eqref{phi}, and $(\omega,c)\in\Omega$, then
\begin{itemize}
  \item[(1)] when $\sigma\in (0,2)$,
 $$
  \|\phi_{\omega,c}\|_{H^1(\R)}\to 0,\quad \mbox{when } c\to -2\sqrt{\omega};
 $$
  \item[(2)] when $\sigma\ge 2$, there exists a constant $c_0=c_0(\sigma)>0$, such that for any $(\omega,c)\in\Omega$,
$$
\|\phi_{\omega,c}\|_{\dot H^{s_c}(\R)}\ge c_0.
$$
\end{itemize}
\end{thm}
%Shown in Appendix, we prove that when $0<\sigma<2$,
%$$
%\|\phi_{\omega,c}\|_{H^1(\R)}\to 0,\quad \mbox{when } c\to -2\sqrt{\omega}.
%$$
%However, when $\sigma\ge 2$, there exists a positive constant $c_0$, such that
%$$
%\|\phi_{\omega,c}\|_{L^2(\R)}\ge c_0,\quad \mbox{for any } (\omega,c)\in \big\{(\omega,c): c^2<4\omega, \mbox{or } 0<c=2\sqrt{\omega}\big\}.
%$$
Hence, the small data scattering is not true for all $0<\sigma<2$, but from (2) in Theorem \ref{thm:main1}, it is reasonable to conjecture that the small data scattering holds when $\sigma\ge 2$. Our second result aims to show this assertion. Note that we can not replace the norm $\dot H^{s_c}$ by $L^{p_c}, p_c=2\sigma$ in Theorem \ref{thm:main1} (2), although the latter space is also invariant under the scaling \eqref{eqs:scaling-p} (see Remark \ref{rem:3.1} below).

Before stating our second main result, we define the working space
\begin{align}
\|u\|_{X_T}=& \|u\|_{L_t^\infty H_x^s([0,T]\times \R)}
+ \|\partial_xu\|_{L_x^\infty L_t^2(\R\times [0,T] )}
+ \sup_{q\in [4, N_0]}\|u\|_{L_x^q L_t^\infty(\R\times [0,T] )}\nonumber\\
&+\|u\|_{L_t^4 L_x^\infty([0,T]\times \R)}
+ \big\|D^{s-\frac 12}u\big\|_{L_x^4 L_t^\infty(\R\times [0,T] )}\nonumber\\
&+\big\|D^{s-\frac 12}\partial_xu\big\|_{L_x^\infty L_t^2(\R\times [0,T] )}
+ \big\|D^{s-\frac 12}u\big\|_{L_t^4 L_x^\infty([0,T]\times \R)}.\label{19:32}
\end{align}
Here $N_0$  is any fixed arbitrary large parameter.
Now our second main result is
\begin{thm}\label{Thm}
Let $\sigma\ge 2$, $\frac 12\leq s \leq1$ and $\varphi \in H^s(\R)$.
Then there exists a constant $\delta_0>0$, such that if $\|\varphi\|_{H^s(\R)}\leq \delta_0$, then the corresponding solution $u$ is global, and
\begin{align*}
\|u\|_{X_\infty}\lesssim  \|\varphi\|_{H^s(\R)}.
\end{align*}
Moreover,  there exists a unique $u_{\pm}$ such that for any $0\leq s'< s$,
\begin{align*}
\big\|u(t)-e^{it\Delta}u_{\pm}\big\|_{H^{s'}(\R)} \to 0 \qquad as \quad t \to \pm \infty.
\end{align*}
\end{thm}

\begin{remark}
The same result is also true when we consider the nonlinearity $P(u,\bar u, \partial_x u, \partial_x\bar u)$ and $d\ge 5$ in \eqref{14:27}
with $P$ has only one derivative. As a comparable result, Pornnopparath \cite{Pornnopparath-PRE} proved that when $\sigma \ge 2$, and is an integer, the equation
in \eqref{GDNLS} is almost globally well-posed in $H^s(\R),s>\frac 12$. Here ``almost" is in the sense that given an arbitrary large $T>0$,
there exists a constant $C=C(T)>0$, such that for any initial data $u_0:\|u_0\|_{H^s(\R)}\leq C$, the corresponding solution is in $[0,T]$.
Theorem \ref{Thm} improves Pornnopparath's result. On one hand, we do not restrict that $\sigma$ is an integer. On the other hand, as a byproduct of scattering, we prove the global well-posedness in $H^s(\R),\frac 12\leq s\leq1$, which contains the ``endpoint" case $s=\frac 12$ and the global well-posedness in the general sense.

We believe that the index $s=\frac 12$ is optimal for local well-posedness in the sense of  uniform continuity of the solution flow.   However, it is not proved in this paper and leaves us an interesting problem to pursue later.

Moreover, it is worth noting that our scattering result  is not applicable in $H^s$ if the initial data $\varphi\in H^s(\mathbb R)$
since the index $s'$ in Theorem \ref{Thm} satisfies $ s'<s$.
\end{remark}

%in which Santos proved the existence and uniqueness of solutions $u\in C([0,T]; H^{\frac 12}(\mathbb R))$ for sufficient small initial data in the case of $\sigma>1$ for \eqref{GDNLS},

%
%\begin{remark}
%The motivation to prove scattering is a serve of our further study. In order to consider the long-time behavior of the solution
%to (gDNLS), the small data scattering theory is initially needed, for example, long-time perturbation theory when we use the
%Concentration-Compactness argument.
%%We believe that the scattering
%%statement of (gDNLS) in $H^s(\R),s\ge 0$ may be not true when $\sigma < 2$. Indeed, as an comparable model, if we consider the power-type nonlinear Schr\"odinger equation
%%$$
%%i\partial_t u+\partial_{xx} u + |u|^{2\sigma} u=0,
%%$$
%%the scattering
%%statement is not correct in $H^s(\R),s\ge 0$ when $\sigma<2$. However, it is not proved in this paper and leaves us an interesting problem to pursue later.
%\end{remark}

Based on the local well-posedness result of Santos \cite{Santos-JDE-2015}, we use the  bootstrap argument to prove Theorem \ref{Thm}. More precisely, by defining the working space $X_T$ as above with any fixed time $T$, our purpose is to show the uniform-in-time estimate:
\begin{align}
\|u\|_{X_T}\leq C_1\|\varphi\|_{H^s(\R)}+C_2\|u\|_{X_T}^{2\sigma+1}.  \label{3:53}
\end{align}
Here $C_1,C_2$ are the constants independent of $T$.
%Hence, as another byproduct, we can prove the local well-posendess for large initial data in $H^s(\R), \frac12\leq s\leq1$, by using the standard fixed point argument. This improves the previous work of Hao \cite{Hao-CPAA-07}.
The tools we use in the present paper are the smoothing effects and the maximal function estimates. Compared with the low power case $\sigma<2$, the maximal function estimates in the case of $\sigma\ge 2$ provide many benefits. This enables us to handle the nonlinearity properly and establish the uniform-in-time estimate.
However, since our desired result is stronger than the previous ones, the situation here has more obstacles. The key ingredients in our proofs are presented below.

(1) A suitable working space is constructed. In order to establish the uniform estimation on time $T$, a related complicated working space need to be constructed. We define the working space $\|u\|_{X_T}$ in \eqref{19:32}.
We shall prove that the estimation of each norm in $X_T$ is closed. The selection of norms plays an important role in our paper.

(2) A key split on the terms involved the fractional derivatives is carried out.
The endpoint Kato-Ponce inequality recently proved by Bourgain and Li \cite{BoLi-DIE-14} shall be used to deal with
some $L^\infty$-$L^\infty$ type Leibniz rule for fractional derivatives.
Moreover, a regular process using H\"older's inequality fails to control these terms by $\|u\|_{X_T}$, since most of the mixed norms like $\sup_{q\in [4, N_0]}\|u\|_{L_x^q L_t^\infty(\R\times [0,T] )}$ are the norm of time ahead. So the subtle  split is established, thus we are able to change the order of the mixed norm in some applicable way.
This idea has significant influence to obtain our whole estimation on the form of $\|u\|_{X_T}$.

The rest of the paper is organized as follows.
In Section \ref{sec:notations}, we give some basic notations and some preliminary estimates that will be used throughout in our paper.
In Section 3, we prove non-scattering result for (gDNLS) in $H^1(\R)$ when $\sigma\in (0,2)$.
In Section 4, we prove scattering result for (gDNLS) in $H^s(\R)$ with small initial datum when $\sigma\ge 2$.
\vskip .2cm

\section{Notation and Preliminary}\label{sec:notations}

\vskip 0.2cm
\subsection{Notation}

\quad We write $X \lesssim Y$ or $Y \gtrsim X$ to indicate $X \leq CY$ for
some constant $C>0$. The notation $a+$ denotes $a+\varepsilon $ for any small $\varepsilon$, and also $a-\varepsilon$
for $a-$. Denote $\langle\cdot\rangle=(1+|\cdot|^2)^{\frac 12}$ and $D^\alpha  =(-\partial_x^2)^{\frac \alpha 2} $. The Hilbert space $H^s(\R)$ is a Banach space of elements such that
$\mathscr \langle\xi\rangle^s \hat u  \in L^2(\R)$, where $\mathscr{F}$ denotes the Fourier transform $\mathscr{F}u (\xi)= \hat u(\xi)
= \int_\R e^{-2\pi i x\cdot \xi} u(x) \dx$, and equipped with the norm $\|u\|_{H^s}= \|\langle\xi\rangle^s  \hat u  (\xi)  \|_{L^2}$. We also have an embedding theorem
that $\|u\|_{H^{s_1}}\lesssim \|u\|_{H^{s_2}}$ for any $s_1 \leq s_2$, $s_1, s_2\in \R$. Throughout the whole paper, the letter $C$ will denote various positive constants
which are of no importance in our analysis.
We use the following norms to denote the mixed spaces $L^q_tL^r_x([0,T]\times\R)$ and $L^r_xL^q_t(\R\times[0,T])$, that is,
\begin{align*}
\|u\|_{L^q_t L_x^r([0,T]\times\R)}=\Big(\int_0^T\|u\|_{L_x^r(\R)}^q\dt\Big)^{\frac 1 q}
\end{align*}
and
\begin{align*}
\|u\|_{L^r_x L_t^q(\R\times[0,T])}=\Big(\int_\R \|u\|_{L_t^q([0,T])}^r\dx\Big)^{\frac 1 r}.
\end{align*}

\subsection{Preliminary}\label{sec:lnear}

\quad In this section, we state some preliminary estimates of the linear Schr\"odinger operator $e^{it\Delta}$ which will be used in our later sections. Firstly, we recall the well-known Strichartz estimates.

\begin{lem}\label{lem:strichartz}
(Strichartz's estimates, see \cite{Cazenave--03}). Let $I\subset \R$ be an interval. For all admissible pair $(q_j,r_j),j=1,2,$  satisfying
\begin{align*}
2\leq q_j,r_j\leq \infty\quad and\quad \frac 2 q_j=\frac12-\frac 1r_j,
\end{align*}
the following estimates hold:
\begin{align}\label{Strichartz1}
\|e^{it\Delta}f \|_{L_t^{q_j} L_x^{r_j}(I\times \R)} \lesssim  \| f  \|_{L^2(\R)};
\end{align}
and
\begin{align}\label{Strichartz2}
\Big\|\int_0^t e^{i(t-t')\Delta} F(x,t')\dt'\Big\|_{L_t^{q_1} L_x^{r_1}(I\times \R)}
\lesssim\|F\|_{L_t^{q_2'} L_x^{r_2'}(I\times \R)},
\end{align}
where $\frac 1 {q_2}+\frac 1{q_2'}=\frac 1 {r_2}+\frac 1{r_2'}=1$.
\end{lem}

The next lemma is the smoothing effects.
\begin{lem}
(Smoothing effects, see \cite{KenigPonceVega-CPAM-93, LinaresPonce-09}). Let $I\subset \R$ be an interval, including $I=R$. Then\\
1)
\begin{align}\label{Smooth1}
\|D^{\frac 12}e^{it\Delta}f \|_{L_x^\infty L_t^2(\R\times I )} \lesssim  \| f  \|_{L^2(\R)}
\end{align}
for all $f\in L^2(\R)$; and\\
2)
\begin{align}\label{smooth2}
\Big\|  D^{\frac 12} \int_0^t e^{i(t-t')\Delta} F(x,t')\dt'\Big\|_{L_t^\infty L_x^2(I\times \R)}
\lesssim           \|F\|_{L_x^1 L_t^2(\R\times I)};
\end{align}
3)
\begin{align}\label{smooth3}
\Big\|  \partial_x \int_0^t e^{i(t-t')\Delta} F(x,t')\dt' \Big\|_{L_x^\infty L_t^2(\R\times I)}
\lesssim           \|F\|_{L_x^1 L_t^2(\R\times I)};
\end{align}
for all $F\in L_x^1 L_t^2(\R\times I)$.
\end{lem}

Next, we introduce the following maximal function estimates for the linear Schr\"odinger equation.
\begin{lem}\label{lem:Maximal}
(Maximal function estimates, see \cite{KenigPonceVega-IUMJ-91, KenigRuiz-TAMC-83, MolinetRibaud-JMPA-04, RogersVillarroya-AASFM-07, Santos-JDE-2015}). Let $I\subset \R$ be an interval. Let $4\leq p<\infty$ and $s\ge \frac 12-\frac 1p$. Then we have
\begin{align}\label{Maximal1}
\|  e^{it\Delta}f     \|_{L_x^p L_t^\infty(\R\times I)}    \lesssim    \|f \|_{H^s(\R)};
\end{align}
and
\begin{align}\label{Maximal2}
\Big\|   \int_0^t e^{i(t-t')\Delta} F(x,t')\dt'\Big\|_{L_x^p L_t^\infty(\R\times I)}
\lesssim           \|F\|_{L_x^1 L_t^2(\R\times I)}+\|F\|_{L_t^1 L_x^2(I\times \R)}.
\end{align}
\end{lem}

Next, we show the Leibniz and chain rule for fractional derivatives, see \cite{BoLi-DIE-14, KenigPonceVega-CPAM-93, LiDong-PRE} and the references therein.
\begin{lem}\label{frac-chain}
(Leibniz and chain rule for fractional derivatives).
Let $I\subset \R$ be an interval. Then
\\

\noindent 1) Let $s \in (0,1)$, $1<p\le \infty$, and $1<p_1,p_2,p_3, p_4 \le \infty$ with $\frac1p=\frac1{p_1}+\frac1{p_2}$, $\frac1p=\frac1{p_3}+\frac1{p_4}$, and let $f,g\in \mathcal S(\R)$,  then
\begin{align}\label{1:04}
\big\|D^s(fg)\big\|_{L^p(\R)}\lesssim \big\|D^sf\big\|_{L^{p_1}(\R)}\|g\|_{L^{p_2}(\R)}+ \big\|D^sg\big\|_{L^{p_3}(\R)}\|f\|_{L^{p_4}(\R)}.
\end{align}
\\
2) Let $s \in (0,1)$ and $p,q,p_1,p_2,q_2\in (1,\infty), q_1\in (1,\infty]$ such that
\begin{align*}
\frac 1p=\frac 1{p_1}+\frac 1{p_2}\quad \text{and}\quad \frac 1q=\frac 1{q_1}+\frac 1{q_2}.
\end{align*}
Then
\begin{align}\label{15:45}
\|D^sF(f)\|_{L_x^pL_t^q(\R\times I)}\lesssim \|F'(f)\|_{L_x^{p_1}L_t^{q_1}(\R\times I)} \|D^s f\|_{L_x^{p_2}L_t^{q_2}(\R\times I)}.
\end{align}
\\
3) Let $s \in (0,1),s_1,s_2\in[0,s]$ with $s=s_1+s_2$. Let $p,p_1,p_2,q,q_1,q_2\in (1,\infty)$ be such that
\begin{align*}
\frac 1p=\frac 1{p_1}+\frac 1{p_2}\quad \text{and}\quad \frac 1q=\frac 1{q_1}+\frac 1{q_2}.
\end{align*}
Then
\begin{align}\label{10:03}
\|D^s(fg)-fD^sg-gD^sf\|_{L_x^pL_t^q(\R\times I)}\lesssim \|D^{s_1}f\|_{L_x^{p_1}L_t^{q_1}(\R\times I)} \|D^{s_2}g\|_{L_x^{p_2}L_t^{q_2}(\R\times I)}.
\end{align}
Moreover, for $s_1=0$ the value $q_1=\infty$ is allowed.
\\

\noindent4) Let $s \in (0,1),s_1,s_2\in[0,s]$ with $s=s_1+s_2$. Let $p_1,p_2,q_1,q_2\in (1,\infty)$
with $1=\frac 1{p_1}+\frac 1{p_2}$ and $\frac 12=\frac 1{q_1}+\frac 1{q_2}$. Then
\begin{align}\label{14:41}
\|D^s(fg)-fD^sg-gD^sf\|_{L_x^1L_t^2(\R\times I)}\lesssim \|D^{s_1}f\|_{L_x^{p_1}L_t^{q_1}(\R\times I)} \|D^{s_2}g\|_{L_x^{p_2}L_t^{q_2}(\R\times I)}.
\end{align}
\end{lem}

\vskip 0.2cm

\section{Proof of Theorem \ref{thm:main1}}

In this section, we consider the solitary wave solutions described in Introduction, and give the proof of Theorem \ref{thm:main1}. Let $\Omega=\{(\omega,c): c\in(-2\sqrt \omega, 2\sqrt \omega)\}$.

\begin{proof}
Note that $\phi_{\omega,c}$ is the solution of the following equation
\begin{align*}
-\partial_x^2\phi+\omega \phi+c i\partial_x\phi-i| \phi|^{2\sigma}\partial_x\phi=0.
\end{align*}
Multiplying on both sides with $\overline{x\partial_x\phi_{\omega,c}}$, taking the real part and integrating over $\R$, we obtain that for any $(\omega,c)\in\Omega$,
\begin{align}
\|\partial_x\phi_{\omega, c}\|_{L^2}^2=\omega\|\phi_{\omega, c}\|_{L^2}^2.\label{7.23-0420}
\end{align}
Hence, for the statement (1), we only need to consider $\|\phi_{\omega,c}\|_{L^2(\R)}$.

Now we fix $\omega>0$ and denote $\alpha=\sqrt{4\omega-c^2}$. From  \eqref{phi}, we find that
\begin{align}
\int_\R |\phi_{\omega,c}|^2\dx
=&\int_\R   |\varphi_{\omega,c}(x)|^2    \dx\nonumber\\
=&\int_\R \Big\{\frac{(\sigma+1)(4\omega-c^2)}{2\sqrt\omega \cosh(\sigma\sqrt{4\omega-c^2}\,x)-c}\Big\}^{\frac{1}{\sigma}}         \dx\nonumber\\
=&\Big(\frac{\sigma+1}{2\sqrt\omega}\Big)^{\frac{1}{\sigma}}
 \alpha^{\frac2\sigma}\int_\R \Big(\frac1{\cosh(\sigma\alpha x)-\frac{c}{2\sqrt\omega}}\Big)^\frac1\sigma \dx\nonumber\\
=&\frac{2}{\sigma}\Big(\frac{\sigma+1}{2\sqrt\omega}\Big)^{\frac{1}{\sigma}}
 \alpha^{\frac2\sigma-1}\int_0^\infty \Big(\frac1{\cosh x-\frac{c}{2\sqrt\omega}}\Big)^\frac1\sigma \dx\nonumber\\
=&C_{\omega,\sigma} \alpha^{\frac2\sigma-1}\int_0^\infty \Big(\frac1{\cosh x-\frac{c}{2\sqrt\omega}}\Big)^\frac1\sigma \dx,\label{22.20}
\end{align}
where $C_{\omega,\sigma}=\frac{2}{\sigma}\Big(\frac{\sigma+1}{2\sqrt\omega}\Big)^{\frac{1}{\sigma}}
$.
For convenience, we denote
$$
I(c)=\int_0^\infty \Big(\frac1{\cosh x-\frac{c}{2\sqrt\omega}}\Big)^\frac1\sigma \dx.
$$
Moreover, we denote $c_\sigma$ as
$$
c_\sigma=I(-2\sqrt{\omega})=\int_0^\infty \Big(\frac1{\cosh x+1}\Big)^\frac1\sigma \dx,
$$
which makes sense since the last integral above is finite.
Note that $I(c)$ is an increasing function, thus we have that for any $c: -2\sqrt{\omega}<c\leq 0$,
\begin{align}
c_\sigma\le I(c)\le I(0).\label{bound-Ic}
\end{align}
This combining with \eqref{22.20} yields that
\begin{align*}
\int_\R |\phi_{\omega,c}|^2\dx \le C_{\omega,\sigma}I(0) \alpha^{\frac2\sigma-1}\to 0,\quad \mbox{when } c\to -2\sqrt{\omega}.
\end{align*}
This proves the statement (1).

For the statement (2), we split into two cases: $-2z_0\sqrt\omega\le c<2\sqrt\omega$ and $-2\sqrt\omega<c<-2z_0\sqrt\omega$. Here $z_0\in (0,1)$ is a constant close enough to 1 (one may set $z_0=\frac{99}{100}$).

Case 1:  $-2z_0\sqrt\omega\le c<2\sqrt\omega$.
We denote $p_c=2\sigma$, then by Sobolev's inequality, it reduces to show that
there exists a constant $c_0>0$, such that for any $(\omega,c)\in\Omega$,
$$
\|\phi_{\omega,c}\|_{L^{p_c}(\R)}\ge c_0.
$$
From  \eqref{phi}, we get that
\begin{align}
\int_\R |\phi_{\omega,c}|^{p_c}\dx
=&\int_\R   |\varphi_{\omega,c}(x)|^{2\sigma}    \dx\notag\\
=&\int_\R \frac{(\sigma+1)(4\omega-c^2)}{2\sqrt\omega \cosh(\sigma\sqrt{4\omega-c^2}\,x)-c}   \dx\notag\\
=&\frac{\sigma+1}{2\sqrt\omega}
 \alpha^2\int_\R \frac1{\cosh(\sigma\alpha x)-\frac{c}{2\sqrt\omega}} \dx\notag\\
=&\frac{2(\sigma+1)}{\sigma}\frac{\alpha}{2\sqrt\omega} \int_0^\infty \frac1{\cosh x-\frac{c}{2\sqrt\omega}} \dx.\label{22.20-2}
\end{align}
Denote $\beta=\frac{2\sqrt \omega-c}{2\sqrt\omega}$, then $\beta> 0$.  Hence,
\begin{align*}
 \int_0^\infty \frac1{\cosh x-\frac{c}{2\sqrt\omega}} \dx&= \int_0^\infty \frac1{\cosh x-1+\beta} \dx
\ge \frac1{\beta}.
\end{align*}
Hence, this last inequality combining with \eqref{22.20-2} gives that
\begin{align}
\int_\R |\phi_{\omega,c}|^{p_c}\dx
\ge &\frac{2(\sigma+1)}{\sigma}\frac{\alpha}{2\sqrt\omega} \frac1{\beta}
=\frac{2(\sigma+1)}{\sigma}\sqrt{\frac{2\sqrt \omega+c}{2\sqrt \omega-c} }\ge \frac{2(\sigma+1)}{\sigma}\sqrt{\frac{1-z_0}{1+z_0}}:=c_0^{p_c}.
\end{align}
%then we have $1< \beta\le 2$ and
%\begin{align*}
%\int_0^\infty \frac1{\cosh x-1+\beta} \dx&= \int_0^\infty \frac1{\cosh x-1+\beta} \dx
%\ge 1.
%\end{align*}
%Hence, this last inequality combining with \eqref{22.20-2} gives that
%\begin{align}
%\int_\R |\phi_{\omega,c}|^{p_c}\dx
%\ge &\frac{2(\sigma+1)}{\sigma}\frac{\alpha}{2\sqrt\omega}
%=\frac{2(\sigma+1)}{\sigma}\sqrt{\frac{2\sqrt \omega+c}{2\sqrt \omega-c} }\ge \frac{2(\sigma+1)}{\sigma}.
%\end{align}

Case 2: $-2\sqrt\omega<c<-2z_0\sqrt\omega$.
Let $c=-2z\sqrt\omega$, then $z_0< z< 1$. First, we rewrite $\phi_{\omega,c}, \varphi_{\omega,c}$ in the following forms. Let
$$
h_z(x)=\left(\frac{1}{\cosh(2\sigma x)+z}\right)^{\frac{1}{2\sigma}}.
$$
Since $z>0$,  there exist positive constants $c_{1\sigma}, C_{j\sigma}, j=1,2,3$ which are independent of $z$, such that
\begin{align}
c_{1\sigma}\le \big\|h_z\big\|_{L^2}\le C_{1\sigma},\quad
\big\|h_z\big\|_{L^{4\sigma+2}}\le C_{2\sigma},\quad
\big\|\partial_x h_z\big\|_{L^2}\le C_{3\sigma}.\label{bound-hz}
\end{align}
Moreover, we rewrite
\begin{align*}
\varphi_{\omega,c }(x)=\big[2(\sigma+1)\big]^{\frac1{2\sigma}}\omega^{\frac1{4\sigma}}(1-z^2)^{\frac1{2\sigma}}
h_z\big(\sqrt\omega\sqrt{1-z^2}x\big),
\end{align*}
and thus
\begin{align*}
\phi_{\omega,c}(x)=&\big[2(\sigma+1)\big]^{\frac1{2\sigma}}\omega^{\frac1{4\sigma}}(1-z^2)^{\frac1{2\sigma}}
h_z\big(\sqrt\omega\sqrt{1-z^2}x\big)\\
&\quad\cdot\exp{\left\{ - iz\sqrt\omega x-i\sqrt{1-z^2}\int_{-\infty}^{\sqrt\omega\sqrt{1-z^2} x}h_z^{2\sigma}(y)dy\right\}}.
\end{align*}
Denote
$$
g_z(x)=(1-z^2)^{\frac1{2\sigma}}
h_z\big(\sqrt{1-z^2}x\big)\exp{\left\{-i\sqrt{1-z^2}\int_{-\infty}^{\sqrt{1-z^2} x}h_z^{2\sigma}(y)dy\right\}},
$$
then
\begin{align*}
\phi_{\omega,c}(x)=&\big[2(\sigma+1)\big]^{\frac1{2\sigma}}\omega^{\frac1{4\sigma}}\exp\{ - iz\sqrt\omega x\}g_z(\sqrt\omega x).
\end{align*}
Hence, by scaling,  we get that
\begin{align}
\big\|\phi_{\omega,c}\big\|_{\dot H^{s_c}}=\big[2(\sigma+1)\big]^{\frac1{2\sigma}}
\big\|\exp\{ - iz x\}g_z\big\|_{\dot H^{s_c}}.\label{phi-g-scaling}
\end{align}
So it reduces to estimate $\big\|\exp\{ - iz x\}g_z\big\|_{\dot H^{s_c}}$. For this, we have
\begin{align*}
\big\|\exp\{ - iz x\}g_z\big\|_{\dot H^{s_c}}^2
=&\int_\R |\xi|^{2s_c} \big|\hat g_z(\xi+z)\big|^2\,d\xi\\
=&\int_\R |\xi-z|^{2s_c} \big|\hat g_z(\xi)\big|^2\,d\xi\\
\ge &\int_{\{|\xi|\le A_0\sqrt{1-z^2}\}} |\xi-z|^{2s_c} \big|\hat g_z(\xi)\big|^2\,d\xi,
\end{align*}
where $A_0$ is a big constant decided later. Since $0<1-z^2\ll 1$, we further get
\begin{align}
\big\|\exp\{ - iz x\}g_z\big\|_{\dot H^{s_c}}^2
\ge \frac12 &\int_{\{|\xi|\le A_0\sqrt{1-z^2}\}} \big|\hat g_z(\xi)\big|^2\,d\xi.\label{8.17-0421}
\end{align}
Now we claim that by choosing $A_0$ large enough,
\begin{align}
\int_{\{|\xi|\le A_0\sqrt{1-z^2}\}} \big|\hat g_z(\xi)\big|^2\,d\xi
\ge \frac12 c_{1\sigma}^2(1-z^2)^{\frac1{\sigma}-\frac12}.    \label{8.06-0421}
\end{align}
Indeed, on one hand,
\begin{align}
\|g_z\|_{L^2}=&\big\|(1-z^2)^{\frac1{2\sigma}}
h_z\big(\sqrt{1-z^2}x\big)\big\|_{L^2}\notag\\
=&(1-z^2)^{\frac1{2\sigma}-\frac14}\|h_z\|_{L^2}\ge c_{1\sigma}(1-z^2)^{\frac1{2\sigma}-\frac14}.\label{7.39-0421}
\end{align}
On the other hand,
\begin{align*}
\int_{\{|\xi|\ge A_0\sqrt{1-z^2}\}} \big|\hat g_z(\xi)\big|^2\,d\xi
\le A_0^{-2}(1-z^2)^{-1} \big\|\partial_x g_z\big\|_{L^2}^2.
\end{align*}
Moreover,  by \eqref{bound-hz}, there exists $C_{4\sigma}>0$ such that
\begin{align*}
 \big\|\partial_x g_z\big\|_{L^2}
 \le & (1-z^2)^{\frac1{2\sigma}+1}\big\|h_z^{2\sigma+1}(\sqrt{1-z^2}x)\big\|_{L^2}+(1-z^2)^{\frac1{2\sigma}+\frac12}\big\|\partial_xh_z(\sqrt{1-z^2}x)\big\|_{L^2}\\
 \le &(1-z^2)^{\frac1{2\sigma}+\frac34}\big\|h_z\big\|_{L^{4\sigma+2}}^{2\sigma+1}+(1-z^2)^{\frac1{2\sigma}+\frac14}\big\|\partial_xh_z\big\|_{L^2}\\
 \le & C_{4\sigma}(1-z^2)^{\frac1{2\sigma}+\frac14}.
\end{align*}
Hence,
\begin{align*}
\int_{\{|\xi|\ge A_0\sqrt{1-z^2}\}} \big|\hat g_z(\xi)\big|^2\,d\xi
\le C_{4\sigma}^2A_0^{-2}(1-z^2)^{\frac1\sigma-\frac12}.
\end{align*}
Choosing $A_0=\sqrt2c_{1\sigma}^{-1}C_{4\sigma}$, then the last estimate above combining with  \eqref{7.39-0421} gives the claim \eqref{8.06-0421}. Thus combining with \eqref{8.06-0421} and \eqref{8.17-0421}, we get
\begin{align*}
\big\|\exp\{ - iz x\}g_z\big\|_{\dot H^{s_c}}
\ge \frac12 c_{1\sigma}(1-z^2)^{\frac1{2\sigma}-\frac14}.
\end{align*}
Now together with the last estimates above and \eqref{phi-g-scaling}, and noting that $\sigma\ge 2$ and $1-z^2\le 1-z_0^2$, we obtain that
\begin{align}
\big\|\phi_{\omega,c}\big\|_{\dot H^{s_c}}\ge \frac1{2} \big[2(\sigma+1)\big]^{\frac1{2\sigma}}c_{1\sigma}(1-z_0^2)^{\frac1{2\sigma}-\frac14}.
\end{align}
Therefore, we establish the desired result in the second case.
This proves the Theorem \ref{thm:main1}.
\end{proof}

\begin{remark}\label{rem:3.1}
One may find from the computation above that for any $(\omega, c)\in \Omega$ and $c<0$, there exists some $c_\sigma>0$ such that
\begin{align*}
\int_\R |\phi_{\omega,c}|^{p_c}\dx
\le &c_\sigma\frac{\alpha}{2\sqrt\omega}.
\end{align*}
Hence, for any $\sigma>0$,
\begin{align*}
\int_\R |\phi_{\omega,c}|^{p_c}\dx
\to &0, \quad \mbox{as } c\to -2\sqrt\omega.
\end{align*}
Therefore, we can not replace $\dot H^{s_c}$ norm by $L^{p_c}$ norm in Theorem \ref{thm:main1} (2). This may be helpful to understand the structure of the equation.
\end{remark}

\vskip 0.2cm

\section{Proof of Theorem \ref{Thm}}\label{Main}

\vskip 0.2cm

In this section, we give the proof of Theorem \ref{Thm}.
Given  $s\ge \frac 12$ and $\varphi \in H^s(\R)$.
Recall the locally well-posed result of Santos \cite{Santos-JDE-2015}, that is, $u\in C([0,T]; H^{\frac 12}(\mathbb R))$ for sufficient small initial data in the case of $\sigma>1$ for \eqref{GDNLS}.
Based on this, fixing $T>0$, we only need to show the uniform-on-time estimate \eqref{3:53}.
Then the bootstrap argument yields that there exists $\delta_0>0$, such that when $\|\varphi\|_{H^s(\R)}\leq \delta_0$,
\begin{align*}
\|u\|_{X_T}\lesssim \|\varphi\|_{H^s(\R)}
\end{align*}
for any $T\in \R$. In the following, we only consider the positive time. Since the negative time direction can be obtained in the same way.

%{\color{red}Throughout the whole section, the main difficulty lies in the assignment of the suitable index with time and space. Although some of the index look nonstandard, they are used to serve the only aiming to get the closed estimate \eqref{3:53}.}

To show \eqref{3:53}, according to the definition of $\|u\|_{X_T}$, we control the norms in the right--hand side of \eqref{19:32} one by one.

\subsection{Estimates on $\|u\|_{L_t^\infty H_x^s([0,T]\times \R)}$}

In this subsection, we give a priori estimate of the solution in $H^s$, which is important for global well-posedness. Moreover, one may find that its proof also plays a crucial role in the proof of scattering in the end of the section.  Before stepping into the complicated details, we give some remark here.  It is worth noting that the maximal function estimate in Lemma \ref{lem:Maximal} has both local ($p<4$) and global results ($p\ge 4$), see \cite{RogersVillarroya-AASFM-07}. The local version was heavily relied on in previous papers to establish the local well-posedness, see for examples \cite{Pornnopparath-PRE, Santos-JDE-2015}. Unfortunately, we emphasize that to obtain the global well-posedness, the local version can not be used in our estimate.

The main result in this subsection is
\begin{align}\label{13:14}
\|u\|_{L_t^\infty H_x^s([0,T]\times \R)}
\lesssim \|\varphi\|_{H_x^s(\R)}+\|u\|_{X_T}^{2\sigma+1}.
\end{align}
We prove \eqref{13:14} by the following two steps.
\\

\emph{Step 1,\qquad $\|u\|_{L_t^\infty L_x^2([0,T]\times \R)}
\lesssim \| \varphi   \|_{H_x^s(\R)}+ \|u\|_{X_T}^{2\sigma+1}$.}

Using the Duhamel formula
\begin{align}\label{Duhamel}
u(t)= e^{it\Delta}\varphi- \int_0^t e^{i(t-t')\Delta} \big(|u|^{2\sigma}\partial_xu\big)(t') \dt',
\end{align}
 and the Strichartz estimates \eqref{Strichartz1} and \eqref{Strichartz2}, we get
\begin{align}
\|  u   \|_{L_t^\infty L_x^2([0,T]\times \R)}
&\lesssim  \big\|  e^{it\Delta}\varphi   \big\|_{L_t^\infty L_x^2([0,T]\times \R)}+\Big\|   \int_0^t e^{i(t-t')\Delta} \big(|u|^{2\sigma}\partial_x u\big)(t')\dt'\Big\|_{L_t^\infty L_x^2([0,T]\times \R)}\nonumber\\
&\lesssim \| \varphi   \|_{L_x^2(\R)}+ \Big\|  |u|^{2\sigma}\partial_x u  \Big\|_{L_t^1 L_x^2([0,T]\times \R)}.\label{15:12}
\end{align}
Next we consider the term $\Big\|  |u|^{2\sigma}\partial_x u  \Big\|_{L_t^1 L_x^2([0,T]\times \R)}$. We claim that
\begin{align}\label{15:08}
\Big\|  |u|^{2\sigma}\partial_x u  \Big\|_{L_t^1 L_x^2([0,T]\times \R)}\lesssim \|u\|_{X_T}^{2\sigma+1}.
\end{align}
Now we write
\begin{align}
\Big\|  |u|^{2\sigma}\partial_x u  \Big\|_{L_t^1 L_x^2([0,T]\times \R)}
=\Big\|  |u|^2\cdot |u|^{2\sigma-2}\partial_x u  \Big\|_{L_t^1 L_x^2([0,T]\times \R)}.\label{key}
\end{align}
We consider the inner integration $L^2_x$ first. By H\"older's inequality, we have
\begin{align*}
\Big\|  |u|^2\cdot |u|^{2\sigma-2}\partial_x u  \Big\|_{ L_x^2(\R)}
\lesssim  \|u\|_{L_x^\infty(\R)}^2\cdot \Big\||u|^{2\sigma-2}\partial_x u  \Big\|_{ L_x^2(\R)}.
\end{align*}
Hence,
\begin{align}\label{16:59}
\Big\|  |u|^{2\sigma}\partial_x u  \Big\|_{L_t^1 L_x^2([0,T]\times \R)}
\lesssim & \Big\|\|u\|_{L_x^\infty(\R)}^2\cdot \big\||u|^{2\sigma-2}\partial_x u  \big\|_{ L_x^2(\R)}\Big\|_{L_t^1([0,T])}\nonumber\\
\lesssim & \|u\|_{L_t^4 L_x^\infty([0,T]\times \R)}^2 \cdot \Big\||u|^{2\sigma-2}\partial_x u  \Big\|_{L_x^2L_t^2(\R\times [0,T] )}.
\end{align}
For the term $\Big\||u|^{2\sigma-2}\partial_x u  \Big\|_{L_x^2L_t^2(\R \times [0,T])}$, note that $2(2\sigma-2) \ge 4$, by H\"older's inequality again we obtain
\begin{align}
\Big\||u|^{2\sigma-2}\partial_x u  \Big\|_{L_x^2L_t^2(\R \times [0,T])}
&\lesssim  \|u\|_{L_x^{2(2\sigma-2)}L_t^\infty(\R\times[0,T])}^{2\sigma-2}\cdot\|\partial_x u\|_{L_x^\infty L_t^2(\R\times[0,T])}\nonumber\\
&\lesssim \|u\|_{X_T}^{2\sigma-1}.\label{4:47}
\end{align}
Putting this result into \eqref{16:59}, we get
\begin{align}
\Big\|  |u|^{2\sigma}\partial_x u  \Big\|_{L_t^1 L_x^2([0,T]\times \R)}
\lesssim& \|u\|_{L_t^4 L_x^\infty([0,T]\times \R)}^2 \cdot  \|u\|_{X_T}^{2\sigma-1}\label{4:10}\\
\lesssim& \|u\|_{X_T}^{2\sigma+1}.\nonumber
\end{align}
Thus we have proved claim \eqref{15:08}. Then by \eqref{15:12}, we have
\begin{align}\label{18:57}
\|u\|_{L_t^\infty L_x^2([0,T]\times \R)}
\lesssim \| \varphi   \|_{H_x^s(\R)}+ \|u\|_{X_T}^{2\sigma+1}.
\end{align}
Thus we have finished the proof on Step 1.
\\

\emph{Step 2,\qquad $\|  D^s  u   \|_{L_t^\infty L_x^2([0,T]\times \R)}
\lesssim \|\varphi\|_{ H_x^s(\R)}+\|u\|_{X_T}^{2\sigma+1}$.}

Using the Duhamel formula \eqref{Duhamel} and the Strichartz estimate \eqref{Strichartz1}, the smoothing effect \eqref{smooth2}, we have
\begin{align}
\|  D^s  u   \|_{L_t^\infty L_x^2([0,T]\times \R)}
&\lesssim  \big\|  e^{it\Delta} D^s \varphi   \big\|_{L_t^\infty L_x^2([0,T]\times \R)}+\Big\| D^\frac 12  \int_0^t e^{i(t-t')\Delta} D^{s-\frac 12} \big(|u|^{2\sigma}\partial_x u\big)\dt'\Big\|_{L_t^\infty L_x^2([0,T]\times \R)}\nonumber\\
&\lesssim \|   D^s \varphi   \|_{ L_x^2(\R)}
+\Big\| D^{s-\frac 12} \big(|u|^{2\sigma}\partial_x u\big)  \Big\|_{L_x^1 L_t^2(\R \times [0,T])}\nonumber\\
&\lesssim  \|   \varphi   \|_{ H_x^s(\R)}
+\Big\|  D^{s-\frac 12} \big(|u|^{2\sigma}\partial_x u\big)  \Big\|_{L_x^1 L_t^2(\R\times[0,T])}.\label{19:18}
\end{align}
Next we claim that
\begin{align}\label{15:29}
\Big\|  D^{s-\frac 12} \big(|u|^{2\sigma}\partial_x u\big)  \Big\|_{L_x^1 L_t^2(\R\times[0,T])}
\lesssim \|u\|_{X_T}^{2\sigma+1}.
\end{align}
To prove this claim, we split it into two cases: $s=\frac 12$ and $\frac 12<s\leq1$.

Case 1: $s=\frac 12$.

By H\"older's inequality, note that $2\sigma \geq4$ , we have
\begin{align}
\Big\|  |u|^{2\sigma}\partial_x u  \Big\|_{L_x^1 L_t^2(\R\times[0,T])}
\lesssim & \Big\|  \|u\|_{L_t^\infty([0,T])}^{2\sigma} \cdot \|\partial_x u\|_{L_t^2([0,T])}  \Big\|_{L_x^1(\R)}\nonumber\\
\lesssim &  \|u\|_{L_x^{2\sigma}L_t^\infty(\R\times[0,T])}^{2\sigma} \cdot \|\partial_x u\|_{L_x^\infty L_t^2(\R\times[0,T])}\nonumber\\
\lesssim & \|u\|_{X_T}^{2\sigma+1}.\label{4.14}
\end{align}

Case 2:  $\frac 12<s\leq1$.

By the Leibniz rule for fractional derivative \eqref{14:41}, we get
\begin{align}
&\Big\|  D^{s-\frac 12} \big(|u|^{2\sigma}\partial_x u\big)  \Big\|_{L_x^1 L_t^2(\R\times[0,T])}\nonumber\\
\lesssim& \Big\|  D^{s-\frac 12} \big(|u|^{2\sigma}\big) \cdot \partial_x u  \Big\|_{L_x^1 L_t^2(\R\times[0,T])}
+\Big\|   |u|^{2\sigma} \cdot D^{s-\frac 12}\partial_x u  \Big\|_{L_x^1 L_t^2(\R\times[0,T])}\nonumber\\
&\hspace{1cm}+\Big\|   D^{s-\frac 12}\big(|u|^{2\sigma}\big)\Big\|_{L_x^{1+} L_t^{\infty-}(\R\times[0,T])}
\cdot \| \partial_x u  \|_{L_x^{\infty-} L_t^{2+}(\R\times[0,T])}.\label{2:05}
\end{align}
We estimate on terms above one by one.

For the first term $\Big\|  D^{s-\frac 12} \big(|u|^{2\sigma}\big) \cdot \partial_x u  \Big\|_{L_x^1 L_t^2(\R\times[0,T])}$ in \eqref{2:05}, by H\"older's inequality, we have
\begin{align}
&\Big\|  D^{s-\frac 12} \big(|u|^{2\sigma}\big) \cdot \partial_x u  \Big\|_{L_x^1 L_t^2(\R\times[0,T])}\nonumber\\
\lesssim & \Big\|\big\|  D^{s-\frac 12} \big(|u|^{2\sigma}\big)   \big\|_{ L_t^{\infty-}([0,T])} \cdot \|\partial_x u\|_{ L_t^{2+}([0,T])}\Big\|_{L_x^1(\R)}\nonumber\\
\lesssim & \Big\|  D^{s-\frac 12} \big(|u|^{2\sigma}\big)   \Big\|_{ L_x^{1+} L_t^{\infty-}(\R\times[0,T])}
\cdot \|\partial_x u\|_{ L_x^{\infty-} L_t^{2+}(\R\times[0,T])}.\label{21:26}
\end{align}
To the term $\Big\|  D^{s-\frac 12} \big(|u|^{2\sigma}\big)   \Big\|_{ L_x^{1+} L_t^{\infty-}(\R\times[0,T])}$ in \eqref{21:26}, using  \eqref{15:45}, note that $\frac 43 (2\sigma-1)\geq4$ ,we get
\begin{align*}
\Big\|  D^{s-\frac 12} \big(|u|^{2\sigma}\big)   \Big\|_{ L_x^{1+} L_t^{\infty-}(\R\times[0,T])}
\lesssim& \Big\| |u|^{2\sigma-1}\Big\|_{ L_x^\frac 43 L_t^\infty(\R\times[0,T])}\cdot   \Big\| D^{s-\frac 12} u\Big\|_{ L_x^{4+} L_t^{\infty-}(\R\times[0,T])}\\
\lesssim & \|u\|_{X_T}^{2\sigma-1} \cdot \Big\| D^{s-\frac 12} u\Big\|_{ L_x^{4+} L_t^{\infty-}(\R\times[0,T])} .
\end{align*}
By interpolating between $\|D^{s-\frac 12}\partial_x u\|_{ L_x^\infty L_t^2(\R\times[0,T])}$ and $\|u\|_{L_x^4 L_t^\infty(\R\times[0,T])}$, we have
that for some $\theta_1 \in (0,1)$,
\begin{align}
\Big\| D^{s-\frac 12} u\Big\|_{ L_x^{4+} L_t^{\infty-}(\R\times[0,T])}
&\lesssim  \|D^{s-\frac 12}\partial_xu\|_{ L_x^\infty L_t^2(\R\times[0,T])}^{\theta_1}
\cdot \|u\|_{L_x^4 L_t^\infty(\R\times[0,T])}^{1-\theta_1}\nonumber\\
&\lesssim  \|u\|_{X_T}.\label{21:17}
\end{align}
Then
\begin{align}
\Big\|  D^{s-\frac 12} \big(|u|^{2\sigma}\big)   \Big\|_{ L_x^{1+} L_t^{\infty-}(\R\times[0,T])}
\lesssim \|u\|_{X_T}^{2\sigma}.\label{16:40}
\end{align}
To the term $\|\partial_x u\|_{ L_x^{\infty-} L_t^{2+}(\R\times[0,T])}$ in \eqref{21:26}, it follows from the interpolation between\\
$\big\|D^{s-\frac 12}\partial_x u\big\|_{L_x^\infty L_t^2(\R\times[0,T])}$ and
$\big\|D^{s-\frac 12} u\big\|_{L_x^{4} L_t^\infty(\R\times[0,T])}$, that is, there exists $\theta_2\in (0,1)$,
\begin{align}
\|\partial_x u\|_{ L_x^{\infty-} L_t^{2+}(\R\times[0,T])}
\lesssim & \big\|D^{s-\frac 12}\partial_x u\big\|_{L_x^\infty L_t^2(\R\times[0,T])}^{\theta_2}
\cdot \big\|D^{s-\frac 12} u\big\|_{L_x^4 L_t^\infty(\R\times[0,T])}^{1-\theta_2}\nonumber\\
\lesssim & \|u\|_{X_T}.\label{16:41}
\end{align}
Inserting \eqref{16:40} and \eqref{16:41} into \eqref{21:26}, we have
\begin{align}\label{16:36}
\Big\|  D^{s-\frac 12} \big(|u|^{2\sigma}\big) \cdot \partial_x u  \Big\|_{L_x^1 L_t^2(\R\times[0,T])}
\lesssim \|u\|_{X_T}^{2\sigma+1}.
\end{align}
Thus we complete the estimate on the first term of \eqref{2:05}.

For the second term $\Big\|   |u|^{2\sigma} \cdot D^{s-\frac 12}\partial_x u  \Big\|_{L_x^1 L_t^2(\R\times[0,T])}$ in \eqref{2:05}, by H\"older's inequality, note that $2\sigma\geq4$,  we have
\begin{align}
\Big\|   |u|^{2\sigma} \cdot D^{s-\frac 12}\partial_x u  \Big\|_{L_x^1 L_t^2(\R\times[0,T])}
\lesssim & \Big\| \| u \|_{ L_t^\infty([0,T])}^{2\sigma} \cdot  \big\| D^{s-\frac 12}\partial_x u \big\|_{ L_t^2([0,T])} \Big\|_{L_x^1(\R)}\nonumber\\
\lesssim & \| u \|_{ L_x^{2\sigma} L_t^\infty(\R\times[0,T])}^{2\sigma} \cdot \Big\| D^{s-\frac 12}\partial_x u \Big\|_{L_x^\infty L_t^2(\R\times[0,T])}\nonumber\\
\lesssim& \|u\|_{X_T}^{2\sigma+1}.\label{16:37}
\end{align}
Thus the estimate on the second term of \eqref{2:05} is also completed.

For the third term $\Big\|   D^{s-\frac 12}\big(|u|^{2\sigma}\big)\Big\|_{L_x^{1+} L_t^{\infty-}(\R\times[0,T])}
\cdot \| \partial_x u  \|_{L_x^{\infty-} L_t^{2+}(\R\times[0,T])}$ in \eqref{2:05}, using \eqref{16:40} and \eqref{16:41}, we have
\begin{align}\label{2:02}
\Big\|   D^{s-\frac 12}\big(|u|^{2\sigma}\big)\Big\|_{L_x^{1+} L_t^{\infty-}(\R\times[0,T])}
\cdot \| \partial_x u  \|_{L_x^{\infty-} L_t^{2+}(\R\times[0,T])}
\lesssim \|u\|_{X_T}^{2\sigma+1}.
\end{align}
Inserting \eqref{16:36}, \eqref{16:37} and \eqref{2:02} into \eqref{2:05}, we have
\begin{align}\label{2:09}
\Big\|  D^{s-\frac 12} \big(|u|^{2\sigma}\partial_x u\big)  \Big\|_{L_x^1 L_t^2(\R\times[0,T])}
\lesssim \|u\|_{X_T}^{2\sigma+1}.
\end{align}
Owing to the above two cases, we finish the proof of claim \eqref{15:29}.
Putting \eqref{2:09} into \eqref{19:18}, we finish the proof on Step 2.

\subsection{Estimates on $\|\partial_xu\|_{L_x^\infty L_t^2(\R\times[0,T])}$}

\quad Using the Duhamel formula \eqref{Duhamel} and the smoothing effects \eqref{Smooth1} and \eqref{smooth3}, we get
\begin{align*}
\|\partial_xu\|_{L_x^\infty L_t^2(\R\times[0,T])}
&\lesssim  \Big\|  e^{it\Delta}  \partial_x\varphi   \Big\|_{L_x^\infty L_t^2(\R\times[0,T])}+\Big\| \partial_x  \int_0^t e^{i(t-t')\Delta} \big(|u|^{2\sigma}\partial_x u\big)(t')\dt'\Big\|_{L_x^\infty L_t^2(\R\times[0,T])}\nonumber\\
&\lesssim \|    \varphi   \|_{ H_x^{\frac 12}(\R)}
+\Big\|  |u|^{2\sigma}\partial_x u  \Big\|_{L_x^1 L_t^2(\R\times[0,T])}\nonumber.
\end{align*}
By \eqref{4.14}, we have
\begin{align}
 \|\partial_xu\|_{L_x^\infty L_t^2(\R\times[0,T])}\lesssim \|   \varphi   \|_{ H_x^s(\R)}+ \|   u   \|_{X_T}^{2\sigma+1}.\label{4.24}
\end{align}

\subsection{Estimates on $\sup_{q\in [4, N_0]}\|u\|_{L_x^q L_t^\infty(\R\times [0,T] )}$}\label{sec:4.3}
%
%
%\quad {\color{red}In this subsection, we use the global maximal function estimates in Lemma \ref{lem:Maximal} to prove our main result. Note that the maximal function estimates have local and global version (see \cite{RogersVillarroya-AASFM-07}). Here we emphasize that the local version can not be used to obtain our global result.}

By Duhamel's formula \eqref{Duhamel} and the maximal function estimates \eqref{Maximal1} and \eqref{Maximal2}, we have
\begin{align*}
\|u\|_{L_x^q L_t^\infty(\R\times[0,T])}
&\lesssim  \|  e^{it\Delta}  \varphi  \|_{L_x^q L_t^\infty(\R\times[0,T])}+\Big\|   \int_0^t e^{i(t-t')\Delta} \big(|u|^{2\sigma}\partial_x u\big)(t')\dt'\Big\|_{L_x^q L_t^\infty(\R\times[0,T])}\\
&\lesssim \|   \varphi   \|_{H_x^s(\R) }
+\Big\|  |u|^{2\sigma}\partial_x u  \Big\|_{L_x^1 L_t^2(\R\times[0,T])}+\Big\|  |u|^{2\sigma}\partial_x u  \Big\|_{L_t^1 L_x^2([0,T]\times\R)},
\end{align*}
where we have used the condition $s\ge\frac 12\geq\frac 12 -\frac 1q$ in Lemma \ref{lem:Maximal}. By \eqref{15:08} and \eqref{4.14}, we obtain
\begin{align*}
\sup_{q\in [4, N_0]}\|u\|_{L_x^q L_t^\infty(\R\times [0,T] )}
\lesssim \|   \varphi   \|_{ H_x^s(\R)}+ \|   u   \|_{X_T}^{2\sigma+1}.
\end{align*}

\subsection{Estimates on $\|u\|_{L_t^4 L_x^\infty([0,T]\times \R)}$}

\quad By Duhamel's formula \eqref{Duhamel} and the Strichartz estimates \eqref{Strichartz1} and \eqref{Strichartz2}, we get
\begin{align*}
\|u\|_{L_t^4 L_x^\infty([0,T]\times \R)}
&\lesssim  \|  e^{it\Delta}  \varphi   \|_{L_t^4 L_x^\infty([0,T]\times \R)}+\Big\|   \int_0^t e^{i(t-t')\Delta} \big(|u|^{2\sigma}\partial_x u\big)(t')\dt'\Big\|_{L_t^4 L_x^\infty([0,T]\times \R)}\\
&\lesssim \|   \varphi   \|_{ L_x^2(\R)}+ \Big\|  |u|^{2\sigma}\partial_x u  \Big\|_{L_t^1 L_x^2([0,T]\times \R)}.
\end{align*}
By \eqref{15:08}, we have
\begin{align*}
\|u\|_{L_t^4 L_x^\infty([0,T]\times \R)}
\lesssim \|   \varphi   \|_{ L_x^2(\R)}+ \|   u   \|_{X_T}^{2\sigma+1}.
\end{align*}

\subsection{Estimates on $\big\|    D^{s-\frac 12}u      \big\|_{L_x^4 L_t^\infty(\R\times [0,T] )}$}\label{sec4.5}
It is worth noting that the endpoint Kato-Ponce inequality in \cite{BoLi-DIE-14} plays a significant role in our estimates in which we meet the Leibniz rule for fractional derivatives in $L^\infty$. Thanks to  this inequality, we are able to deal with the term $ \| D^{s-\frac 12} (|u|^2)\|_{ L_x^\infty(\R)}$.

Using Duhamel's formula \eqref{Duhamel} and the maximal function estimates \eqref{Maximal1} and \eqref{Maximal2}, we get
\begin{align}
\big\|    D^{s-\frac 12}u      \big\|_{L_x^4 L_t^\infty(\R\times[0,T])}
\lesssim & \big\|  e^{it\Delta}  D^{s-\frac 12}\varphi     \big\|_{L_x^4L_t^\infty(\R\times[0,T])}\nonumber\\
&+\Big\|      \int_0^t  e^{i(t-t')\Delta}  D^{s-\frac 12}\big(   |u|^{2\sigma}\partial_xu     \big)(t')\dt'                \Big\|_{L_x^4L_t^\infty(\R\times[0,T])}\nonumber\\
\lesssim& \big\|    D^{s-\frac 12}\varphi     \big\|_{H_x^\frac 12(\R)}
+\Big\| D^{s-\frac 12} \big(|u|^{2\sigma}\partial_x u\big)  \Big\|_{L_x^1 L_t^2(\R\times[0,T])}\nonumber\\
&+\Big\| D^{s-\frac 12} \big(|u|^{2\sigma}\partial_x u\big)  \Big\|_{L_t^1 L_x^2([0,T]\times\R)},\nonumber\\
\lesssim & \| \varphi\|_{H_x^s(\R)}
+\Big\| D^{s-\frac 12} \big(|u|^{2\sigma}\partial_x u\big)  \Big\|_{L_x^1 L_t^2(\R\times[0,T])}\nonumber\\
&+\Big\| D^{s-\frac 12} \big(|u|^{2\sigma}\partial_x u\big)  \Big\|_{L_t^1 L_x^2([0,T]\times\R)}.\label{12:53}
\end{align}

Recall that the term $\Big\| D^{s-\frac 12} \big(|u|^{2\sigma}\partial_x u\big)  \Big\|_{L_x^1 L_t^2(\R \times [0,T])}$ is already estimated in \eqref{15:29}. So we
only need to consider the term $\Big\| D^{s-\frac 12} \big(|u|^{2\sigma}\partial_x u\big)  \Big\|_{L_t^1 L_x^2([0,T]\times \R)}$.
Now we claim that
\begin{align}\label{17:28}
\Big\| D^{s-\frac 12} \big(|u|^{2\sigma}\partial_x u\big)  \Big\|_{L_t^1 L_x^2([0,T]\times \R)}
\lesssim\|u\|_{X_T}^{2\sigma+1}.
\end{align}
Again, we split it into two cases: $s=\frac 12$ and $\frac 12<s\leq1$.

Case 1: $s=\frac 12$.

The term $\Big\|  |u|^{2\sigma}\partial_x u \Big\|_{L_t^1 L_x^2([0,T]\times \R)}$ is already estimated in \eqref{15:08}.

Case 2: $\frac 12<s\leq1$.

Using a similar
treatment as \eqref{key}, we have
\begin{align*}
\Big\| D^{s-\frac 12} \big(|u|^{2\sigma}\partial_x u\big)  \Big\|_{L_t^1 L_x^2([0,T]\times \R)}
=\Big\| D^{s-\frac 12} \big(|u|^2\cdot|u|^{2\sigma-2}\partial_x u\big)  \Big\|_{L_t^1 L_x^2([0,T]\times \R)}.
\end{align*}
Further, using the Leibniz rule for fractional derivative \eqref{1:04}, we have
\begin{align*}
\Big\| D^{s-\frac 12} \big(|u|^2\cdot|u|^{2\sigma-2}\partial_x u\big)  \Big\|_{ L_x^2(\R)}
\lesssim& \Big\| D^{s-\frac 12} (|u|^2)\Big\|_{ L_x^\infty(\R)}  \cdot   \Big\| |u|^{2\sigma-2}\partial_x u  \Big\|_{ L_x^2(\R)}\\
&+\|  u\|_{ L_x^\infty(\R)}^2  \cdot   \Big\| D^{s-\frac 12}\big(|u|^{2\sigma-2}\partial_x u\big)  \Big\|_{ L_x^2(\R)},
\end{align*}
then
\begin{align*}
&\Big\| D^{s-\frac 12} \big(|u|^{2\sigma}\partial_x u\big)  \Big\|_{L_t^1 L_x^2([0,T]\times \R)}\nonumber\\
\lesssim& \Big\| D^{s-\frac 12} (|u|^2)\Big\|_{L_t^2 L_x^\infty([0,T]\times \R)}  \cdot   \Big\| |u|^{2\sigma-2}\partial_x u  \Big\|_{ L_x^2L_t^2(\R\times[0,T])}\nonumber\\
&+\|  u\|_{ L_t^4 L_x^\infty([0,T]\times \R)}^2  \cdot   \Big\| D^{s-\frac 12}\big(|u|^{2\sigma-2}\partial_x u\big)  \Big\|_{ L_x^2L_t^2(\R\times[0,T])}.
\end{align*}
Note that the term $\Big\| |u|^{2\sigma-2}\partial_x u  \Big\|_{ L_x^2L_t^2(\R\times[0,T])}$ has been considered in \eqref{4:47}, so we only need to deal with the terms $\Big\| D^{s-\frac 12} (|u|^2)\Big\|_{L_t^2 L_x^\infty([0,T]\times \R)}$
and $\Big\| D^{s-\frac 12}\big(|u|^{2\sigma-2}\partial_x u\big)  \Big\|_{ L_x^2L_t^2(\R\times[0,T])}$  respectively.

For the term $\Big\| D^{s-\frac 12} (|u|^2)\Big\|_{L_t^2 L_x^\infty([0,T]\times \R)}$, by \eqref{1:04} and the H\"older inequality, we have
\begin{align*}
\Big\| D^{s-\frac 12} (|u|^2)\Big\|_{L_t^2 L_x^\infty([0,T]\times\R)}
\lesssim& \|u\|_{L_t^4L_x^\infty([0,T]\times \R)}  \cdot  \big\| D^{s-\frac 12} u\big\|_{L_t^4 L_x^\infty([0,T]\times \R)}\nonumber\\
\lesssim & \|u\|^2_{X_T}.
\end{align*}

For the term $\Big\| D^{s-\frac 12}\big(|u|^{2\sigma-2}\partial_x u\big)  \Big\|_{ L_x^2L_t^2(\R\times[0,T])}$,
we claim that
\begin{align}
\Big\| D^{s-\frac 12}\big(|u|^{2\sigma-2}\partial_x u\big)  \Big\|_{ L_x^2L_t^2(\R\times[0,T])}
\lesssim \|u\|_{X_T}^{2\sigma-1}.\label{2:32}
\end{align}
Indeed, using the Leibniz rule for fractional derivative \eqref{10:03}, we obtain
\begin{align*}
&\Big\| D^{s-\frac 12}\big(|u|^{2\sigma-2}\partial_x u\big)  \Big\|_{ L_x^2L_t^2(\R\times[0,T])}\nonumber\\
\lesssim& \Big\| D^{s-\frac 12}\big(|u|^{2\sigma-2}\big)\cdot\partial_x u  \Big\|_{ L_x^2L_t^2(\R\times[0,T])}
+\Big\| |u|^{2\sigma-2} \cdot D^{s-\frac 12}\partial_x u  \Big\|_{ L_x^2L_t^2(\R\times[0,T])}\nonumber\\
&+\Big\| D^{s-\frac 12}\big(|u|^{2\sigma-2}\big) \Big\|_{ L_x^{2+} L_t^{\infty-}(\R\times[0,T])} \cdot \|\partial_x u \|_{ L_x^{\infty-} L_t^{2+}(\R\times[0,T])}.
\end{align*}
To the term $\Big\| D^{s-\frac 12}\big(|u|^{2\sigma-2}\big)\cdot\partial_x u  \Big\|_{ L_x^2L_t^2(\R\times[0,T])}$,
using the H\"older inequality, \eqref{15:45} and \eqref{16:41}, note that $4(2\sigma-3)\geq4$, we have
\begin{align}
&\Big\| D^{s-\frac 12}\big(|u|^{2\sigma-2}\big)\cdot\partial_x u  \Big\|_{ L_x^2L_t^2(\R\times[0,T])}\nonumber\\
\lesssim & \Big\|\big\|  D^{s-\frac 12} \big(|u|^{2\sigma-2}\big)   \big\|_{ L_t^{\infty-}([0,T])} \cdot \|\partial_x u\|_{ L_t^{2+}([0,T])}\Big\|_{L_x^2(\R)}\nonumber\\
\lesssim & \Big\|  D^{s-\frac 12} \big(|u|^{2\sigma-2}\big)   \Big\|_{ L_x^{2+} L_t^{\infty-}(\R\times[0,T])}
\cdot \|\partial_x u\|_{ L_x^{\infty-} L_t^{2+}(\R\times[0,T])}\nonumber\\
\lesssim& \Big\| |u|^{2\sigma-3}\Big\|_{ L_x^4 L_t^\infty(\R\times[0,T])}\cdot   \Big\| D^{s-\frac 12} u\Big\|_{ L_x^{4+} L_t^{\infty-}(\R\times[0,T])}\cdot\|u\|_{X_T}\nonumber\\
\lesssim & \|u\|_{X_T}^{2\sigma-2}\cdot   \Big\| D^{s-\frac 12} u\Big\|_{ L_x^{4+} L_t^{\infty-}(\R\times[0,T])}.\label{2:23}
\end{align}
Hence,
by \eqref{21:17}, we get
\begin{align*}
\Big\| D^{s-\frac 12}\big(|u|^{2\sigma-2}\big)\cdot\partial_x u  \Big\|_{ L_x^2L_t^2(\R\times[0,T])}
\lesssim \|u\|_{X_T}^{2\sigma-1}.
\end{align*}
To the term $\Big\| |u|^{2\sigma-2} \cdot D^{s-\frac 12}\partial_x u  \Big\|_{ L_x^2L_t^2(\R\times[0,T])}$,
using H\"older's inequality, note that $2(2\sigma-2)\geq4$, we get
\begin{align*}
\Big\| |u|^{2\sigma-2} \cdot D^{s-\frac 12}\partial_x u  \Big\|_{ L_x^2L_t^2(\R\times[0,T])}
\lesssim & \| u\|_{ L_x^{2(2\sigma-2)}L_t^\infty(\R\times[0,T])}^{2\sigma-2}  \cdot  \Big\| D^{s-\frac 12}\partial_x u  \Big\|_{ L_x^\infty L_t^2(\R\times[0,T])}\nonumber\\
\lesssim& \|u\|_{X_T}^{2\sigma-1}.
\end{align*}
To the term $\Big\| D^{s-\frac 12}\big(|u|^{2\sigma-2}\big) \Big\|_{ L_x^{2+} L_t^{\infty-}(\R\times[0,T])} \cdot \|\partial_x u \|_{ L_x^{\infty-} L_t^{2+}(\R\times[0,T])}$, as the same estimation in \eqref{2:23}, we have
\begin{align*}
\Big\| D^{s-\frac 12}\big(|u|^{2\sigma-2}\big) \Big\|_{ L_x^{2+} L_t^{\infty-}(\R\times[0,T])} \cdot \|\partial_x u \|_{ L_x^{\infty-} L_t^{2+}(\R\times[0,T])}\lesssim \|u\|_{X_T}^{2\sigma-1}.
\end{align*}
Thus we finish the proof of claim \eqref{17:28} and \eqref{2:32} and obtain
\begin{align*}
 \big\|    D^{s-\frac 12}u      \big\|_{L_x^4 L_t^\infty(\R\times [0,T] )}
\lesssim \| \varphi\|_{H_x^s(\R)}+\|u\|^{2\sigma+1}_{X_T}.
\end{align*}

\subsection{Estimates on $\big\|    D^{s-\frac 12}\partial_xu     \big\|_{L_x^\infty L_t^2(\R\times[0,T])}$}

\quad Using Duhamel's formula \eqref{Duhamel} and the smoothing effects \eqref{Smooth1} and \eqref{smooth3}, we get
\begin{align*}
\big\|    D^{s-\frac 12}\partial_xu        \big\|_{L_x^\infty L_t^2(\R\times[0,T])}
\lesssim  &  \big\| e^{it\Delta}    D^{s-\frac 12}\partial_x\varphi      \big\|_{L_x^\infty L_t^2(\R\times[0,T])}\\
&+\Big\|\partial_x\int_0^t  e^{i(t-t')\Delta}      D^{s-\frac 12}\big(|u|^{2\sigma}\partial_xu\big)(t')      \dt' \Big\|_{L_x^\infty L_t^2(\R\times[0,T])}\\
\lesssim & \big\|D^s\varphi \big\|_{L_x^2(\R)}+\Big\|      D^{s-\frac 12}\big(|u|^{2\sigma}\partial_xu\big) \Big\|_{L_x^1 L_t^2(\R\times[0,T])}\\
\lesssim  &\|\varphi\|_{H_x^s(\R)}+\Big\|      D^{s-\frac 12}\big(|u|^{2\sigma}\partial_xu\big) \Big\|_{L_x^1 L_t^2(\R\times[0,T])}.
\end{align*}
Note that we already have the estimation on $\Big\|      D^{s-\frac 12}\big(|u|^{2\sigma}\partial_xu\big) \Big\|_{L_x^1 L_t^2(\R\times[0,T])}$ in \eqref{2:09}.
Then
\begin{align*}
\big\|    D^{s-\frac 12}\partial_xu      \big\|_{L_x^\infty L_t^2(\R\times[0,T])}
\lesssim \|\varphi\|_{H_x^s(\R)}+ \|u\|^{2\sigma+1}_{X_T}.
\end{align*}

\subsection{Estimates on $\big\| D^{s-\frac 12}u  \big\|_{L_t^4 L_x^\infty([0,T]\times \R)}$}

\quad By Duhamel's formula \eqref{Duhamel} and the Strichartz estimates \eqref{Strichartz1} and \eqref{Strichartz2}, we have
\begin{align*}
\big\| D^{s-\frac 12}u  \big\|_{L_t^4 L_x^\infty([0,T]\times \R)}
\lesssim & \big\|e^{it\Delta} D^{s-\frac 12} \varphi      \big\|_{L_t^4 L_x^\infty([0,T]\times \R)}\\
&+\Big\|\int_0^t  e^{i(t-t')\Delta}     D^{s-\frac 12}\big(|u|^{2\sigma}\partial_xu\big)(t')      \dt' \Big\|_{L_t^4 L_x^\infty([0,T]\times \R)}\\
\lesssim& \big\| D^{s-\frac 12} \varphi      \big\|_{ L_x^2(\R)}
+\big\|    D^{s-\frac 12}\big(|u|^{2\sigma}\partial_xu\big)  \big\|_{L_t^1 L_x^2([0,T]\times \R)}\\
\lesssim& \| \varphi \|_{ H_x^s(\R)}
+\big\|    D^{s-\frac 12}\big(|u|^{2\sigma}\partial_xu\big)  \big\|_{L_t^1 L_x^2([0,T]\times \R)}.
\end{align*}
Note that the estimation on $\big\|    D^{s-\frac 12}\big(|u|^{2\sigma}\partial_xu\big)  \big\|_{L_t^1 L_x^2([0,T]\times \R)}$ is obtained in \eqref{17:28}. Then we have
\begin{align*}
\big\| D^{s-\frac 12}u  \big\|_{L_t^4 L_x^\infty([0,T]\times \R)}
\lesssim \|\varphi\|_{H_x^s(\R)}+ \|u\|^{2\sigma+1}_{X_T}.
\end{align*}

Finally, all the estimates on $ \|u\|_{X_T}$ are obtained and we have
\begin{align*}
\|u\|_{X_T}\lesssim \|\varphi\|_{H_x^s(\R)}+ \|u\|^{2\sigma+1}_{X_T}
\end{align*}
uniformly on $T$.
Hence we get $\|u\|_{X_\infty}\lesssim \|\varphi\|_{H_x^s(\R)}$, which gives the proof of the global well-posedness.

Next we prove the scattering statement. Set
\begin{align*}
u_{+}=\varphi-\int_0^{+\infty}e^{-it'\Delta}\big(|u|^{2\sigma}\partial_xu \big)\, \mathrm{d}t'.
\end{align*}
%By smoothing effects \eqref{smooth2} and \eqref{15:08}, \eqref{15:29}, we obtain
%\begin{align*}
%\|u_+\|_{H_x^s(\R)}&\lesssim\|\varphi\|_{H_x^s(\R)}+\Big\|\int_0^{+\infty}e^{-it'\Delta}\big(|u|^{2\sigma}\partial_xu\big)\, \mathrm{d}t'\Big\|_{L^2(\R)}\\
%&+\Big\| D^\frac12\int_0^{+\infty}e^{-it'\Delta}D^{s-\frac12}\big(|u|^{2\sigma}\partial_xu \big)\, \mathrm{d}t'\Big\|_{L^2(\R)}\\
%&\lesssim\|\varphi\|_{H_x^s(\R)}+\Big\|  |u|^{2\sigma}\partial_x u \Big\|_{L_t^1 L_x^2([0,+\infty]\times \R)}+\Big\|  D^{s-\frac 12} \big(|u|^{2\sigma}\partial_x u\big)  \Big\|_{L_x^1 L_t^2(\R\times[0,+\infty])}\\
%&\lesssim\|\varphi\|_{H_x^s(\R)}+\|u\|^{2\sigma+1}_{X_\infty}\\
%&\lesssim\|\varphi\|_{H_x^s(\R)}+\|\varphi\|^{2\sigma+1}_{H_x^s(\R)}.
%\end{align*}
%Then we have
%\begin{align*}
%\|u(t)-e^{it\Delta}u_{+}\|_{H_x^{s}(\R)}&\lesssim\|u\|_{H_x^s(\R)}+\|u_+\|_{H_x^s(\R)}\\
%&\lesssim\|u\|_{X_T}+\|\varphi\|_{H_x^s(\R)}+\|\varphi\|^{2\sigma+1}_{H_x^s(\R)}\\
%&\lesssim\|\varphi\|_{H_x^s(\R)}+\|\varphi\|^{2\sigma+1}_{H_x^s(\R)}.
%\end{align*}
Using Duhamel's formula \eqref{Duhamel}, we have
\begin{align*}
u(t)-e^{it\Delta}u_{+}=\int_t^{+\infty}e^{i(t-t')\Delta}\big(|u|^{2\sigma}\partial_xu \big)\, \mathrm{d}t'.
\end{align*}
By interpolation, for any $0\leq s'<s$, we have that for some $\theta \in[0,1)$ ,
\begin{align}
\|u(t)-e^{it\Delta}u_{+}\|_{H^{s'}(\R)}&\lesssim\|u(t)-e^{it\Delta}u_{+}\|_{L^2(\R)}^{\theta}\cdot\|u(t)-e^{it\Delta}u_{+}\|_{\dot H_x^{s}(\R)}^{1-\theta}.\label{8.1}
\end{align}
For the term $\|u(t)-e^{it\Delta}u_{+}\|_{L^2(\R)}$, by \eqref{4:10}, we have
\begin{align*}
\|u(t)-e^{it\Delta}u_{+}\|_{L^2(\R)}\lesssim&\Big\|  |u|^{2\sigma}\partial_x u  \Big\|_{L_t^1 L_x^2([t,+\infty]\times \R)}\nonumber\\
\lesssim&\|u\|_{L_t^4 L_x^\infty([t,+\infty]\times \R)}^2 \cdot  \|u\|_{X_\infty}^{2\sigma-1}\nonumber.
\end{align*}
Since $\|u\|_{X_\infty}\lesssim \|\varphi\|_{H_x^s(\R)}$, we get
\begin{align*}
\|u\|_{L_t^4 L_x^\infty([t,+\infty]\times \R)}\rightarrow0,\quad\quad\quad when \quad t\rightarrow+\infty.
\end{align*}
Therefore
\begin{align}
\|u(t)-e^{it\Delta}u_{+}\|_{L^2(\R)}\rightarrow0,\quad\quad when \quad t\rightarrow+\infty.\label{8.77}
\end{align}
For the term $\|u(t)-e^{it\Delta}u_{+}\|_{\dot H_x^{s}(\R)}$, by smoothing effects \eqref{smooth2} and \eqref{15:29} , we have
\begin{align*}
\|u(t)-e^{it\Delta}u_{+}\|_{\dot H_x^{s}(\R)}&\lesssim\Big\|D^\frac12\int_t^{+\infty}e^{i(t-t')\Delta}D^{s-\frac12}\big(|u|^{2\sigma}\partial_xu \big)\, \mathrm{d}t'\Big\|_{L^2(\R)}\nonumber\\
&\lesssim\Big\|  D^{s-\frac 12} \big(|u|^{2\sigma}\partial_x u\big)  \Big\|_{L_x^1 L_t^2(\R\times[0,+\infty])}\nonumber\\
&\lesssim\|u\|_{X_\infty}^{2\sigma+1}.
\end{align*}
Hence, we have
\begin{align}
\|u(t)-e^{it\Delta}u_{+}\|_{\dot H_x^{s}(\R)}\lesssim\|\varphi\|_{H_x^s(\R)}^{2\sigma+1}.\label{8.3}
\end{align}
Combining the estimates \eqref{8.1}, \eqref{8.77} and \eqref{8.3}, we get
\begin{align*}
\|u(t)-e^{it\Delta}u_{+}\|_{H^{s'}(\R)}\rightarrow0,\quad\quad\quad\quad as\quad t\rightarrow+\infty.
\end{align*}
This proves the scattering statement and thus finish the proof of the Theorem \ref{Thm}.

\section*{Acknowledgements}
The authors are grateful to the referee who gave useful notes on grammatical/typographical errors and many helpful comments and suggestions.

\end{document}